\documentclass[10pt]{article}

\usepackage{amsfonts}
\usepackage{epsfig}
\usepackage{graphicx}
\usepackage{amsmath}
\usepackage{amssymb}

\newcounter{main}

\newtheorem{theorem}{Theorem}[section]
\newtheorem{proposition}[theorem]{Proposition}
\newtheorem{lemma}[theorem]{Lemma}

\newtheorem{maintheorem}{Theorem}

\newcommand{\blanksquare}{\,\,\,$\sqcup\!\!\!\!\sqcap$}
\newenvironment{proof}{{\flushleft {\bf Proof: }}}{\blanksquare}

\newcounter{example}
{{\stepcounter{example}}{\flushleft {\bf Example \arabic{example}:}}}%
{\par}

\title{\textbf{Removing zero Lyapunov exponents in volume-preserving flows}}

\author{ M\'{a}rio Bessa \thanks{Supported by FCT-FSE, SFRH/BPD/20890/2004. }
\space and Jorge Rocha \thanks{Partially supported by FCT-POCTI/MAT/61237/2004.}}
\date{October 17, 2006}

\begin{document}
\maketitle

\begin{abstract}
Baraviera and Bonatti in ~\cite{BB} proved that it is possible to perturb, in the $C^{1}$-topology, a volume-preserving and partial hyperbolic diffeomorphism  in order to obtain a non-zero sum of all the Lyapunov exponents in the central direction. In this article we obtain the analogous  result for volume-preserving flows.
\end{abstract}

\bigskip

\emph{keywords:} Dominated splitting; volume-preserving flows; Lyapunov exponents; Stably ergodic.\\

\begin{section}{Introduction}
One of the central problems in dynamical systems is to understand how orbits of nearby points asymptotically diverge (or converge) as one iterates a linear approximation of a system. These rates are measured by the Lyapunov exponents. Given an invariant measure for a flow the Oseledets's Theorem~\cite{O} assures that Lyapunov exponents exist at least for almost every points relatively to that measure. Moreover, if all these exponents are different from zero and if a  H\"{o}lder regularity condition on the derivative is satisfied then \emph{Pesin Theory} of non-uniformly
hyperbolic systems guarantees a very complete invariant manifold theory for a.e. point in $M$ relatively to that measure. These geometric objects reveal to be crucial in the development of the smooth ergodic theory. Hence it is very important to establish favorable conditions on which we can approximate our original system in order to remove zero Lyapunov exponents.

In an opposite direction (see~\cite{Be}) it is proved that, in $C^1$-topology, zero Lyapunov exponents are abundant, which gives us a reason why given a conservative flow it is so difficult to detect positive Lyapunov exponents.

In~\cite{BB} Baraviera and Bonatti succeeded in removing zero Lyapunov exponents by small $C^{1}$-perturbations in discrete stably ergodic dynamical systems with some weak hyperbolicity. Here we give a contribution to the theory by establishing an analogous result for the continuous-time setting.  

In order to state our main result let us first introduce some notations and definitions. 

Let $M$ be a $n$-dimensional compact, connected, without boundary $C^{\infty}$ Riemannian
manifold. Let $\mu$ be the measure induced by a volume-form defined in $M$. We call $\mu$ the Lebesgue measure. Let $\mathfrak{X}^{1}_{\mu}(M)$ be the set of $C^1$-vector fields $X:M\rightarrow{TM}$ such that $\mu$ is $X^{t}$-invariant, where $X^{t}:M\rightarrow M$ is the flow associated to $X$. We endowed $\mathfrak{X}^{1}_{\mu}(M)$ with the $C^{1}$-norm. We say that $X \in \mathfrak{X}^{1}_{\mu}(M)$ is \emph{stably ergodic} (with respect to the $C^1$-topology) if there exists a $C^1$ neighborhood of $X$, $\mathcal W$, such that $\mu$ is an ergodic invariant measure for every $Y \in \mathcal W$.

We recall that the tangent maps $DX:{TM}\rightarrow T(TM)$ and $DX^{t}:TM\rightarrow TM$ are related: the tangent flow $DX^{t}_{p}$ satisfies the non-autonomous linear variational equation $\dot{u}(t)=DX_{X^{t}(p)}\cdot u(t)$. 

We denote by $Sing(X)$ the set of singularities of $X$, say the points $x\in M$ such that $X(x)=\vec{0}$, and by $R:=M\setminus Sing(X)$ the set of regular points. 
Given $x\in R$ we consider its normal bundle $N_{x}=X(x)^{\perp}\subset T_{x}M$ and define the linear Poincar\'{e} flow by $P_{X}^{t}(x):=\Pi_{X^{t}(x)}\circ DX^{t}_{x}$ where $\Pi_{X^{t}(x)}:T_{X^{t}(x)}M\rightarrow N_{X^{t}(x)}$ is the projection along the direction of $X(X^{t}(x))$. An $P_{X}^{t}$-invariant splitting $N=N^{1}\oplus ...\oplus N^{k}$ is called a $m$-\emph{dominated splitting} for the linear Poincar\'{e} flow if there exists $m \in {\mathbb N}$ such that, for all $x\in M$ and $0\leq i<j \leq k $, we have:
$$\frac{\|P_{X}^{m}(x)|_{N^{j}_{x}}\|}{{\mathfrak m}(P_{X}^{m}(x)|_{N^{i}_{x}})}\leq \frac{1}{2},$$
where ${\mathfrak m}(\cdot)$ denotes the co-norm of an operator, that is ${\mathfrak m}(A)= \| A^{-1} \|^{-1}$. We say that the subbundle $N_{i}$ is \emph{hyperbolic} if there exists $k\in\mathbb{N}$ such that either $\|(P_{X}^{k}(x)\cdot u)^{-1}\|\leq 1/2$ (expanding), for all $x\in M$ and any unit vector $u\in N_{i}(x)$, or $\|P_{X}^{k}(x)\cdot u\|\leq 1/2$ (contracting), for all $x\in M$ and any unit vector $u\in N_{i}(x)$.

Given a vector field $X$, an $X^t$-invariant set $\Lambda \subseteq M\setminus Sing(X)$ is (uniformly)~\emph{partially hyperbolic} if there exists an $P_{X}^{t}$-invariant dominated splitting $N=N^{u}\oplus N^c \oplus N^{s}$ in $\Lambda$ such that  $N^u$ is hyperbolic expanding and $N^s$ is hyperbolic contracting, moreover these two subbundles are not simultaneously trivial. We observe that if $\Lambda$ is compact then this definition is equivalent to the analogous one by considering the dominated splitting $T_\Lambda M=E^{u}\oplus E^c \oplus E^{s}$ invariant by the tangent flow, $DX^t$.

A consequence of the  Oseledets's Theorem is that, given $X\in\mathfrak{X}^{1}_{\mu}(M)$, Lebesgue a.e. point $p\in{M}$ admits a (Oseledets's) splitting $N_{p}=N^{1}_{p}\oplus...\oplus{N^{k(p)}_{p}}$ and Lyapunov exponents $\lambda_{1}(p)>...>\lambda_{k(p)}(p)$ such that $P_{X}^{t}(p)(N^{i}_{p})=N^{i}_{X^{t}(p)}$ and $\underset{t\rightarrow{\pm{\infty}}}{\text{lim}}\frac{1}{t}\log{\|P_{X}^{t}(p)\cdot v^{i}\|=\lambda_{i}(p)}$
for any $v^{i}\in{N^{i}_{p}\setminus\{0\}}$ and $i=1,...,k(p)$ (see for exam\-ple~\cite{Be}). Moreover we also have
\begin{equation}\label{angle}
\lim_{t\rightarrow{\pm{\infty}}}\frac{1}{t}\log{|\text{det}(P_{X}^{t}(p))|=\sum_{i=1}^{k(p)}\lambda_{i}(p).dim(N^{i}_{p})},
\end{equation}
Finally, the sum of all central Lyapunov exponents of $X$ is given by $$\Sigma^c(X)=\int_M \text{log} | \text{det} P_X^1(p) _{|N^c} | d\mu(p).$$

We observe that the ergodicity of the flow implies that the Lyapunov exponents and the dimensions of the associated subbundles are a.e.-constant, and note that the Birkoff's ergodic Theorem and (\ref{angle}) assure that $\Sigma^c(X)=\sum_{i=1}^{k}\lambda_{i}.dim(N^{i})$.

Let us now state the fundamental result of this article.

\begin{maintheorem}\label{teo1}
Let $X \in \mathfrak{X}^{1}_{\mu}(M)$  be a stably ergodic flow and assume that all the singularities of $X$ are linear hyperbolic and that $R=M\setminus Sing(X)$ is a partially hyperbolic set. Then, either $\Sigma^c(X)\not= 0$, or else $X$ may be approximated, in the $C^{1}$-topology, by $Y \in \mathfrak{X}^{2}_{\mu}(M)$ for which $\Sigma^c(Y) \not=0$.
\end{maintheorem}

\end{section}

\begin{section}{Outline of the proof  and a perturbation scheme}

Let us give a brief sketch of the proof. Fix a vector field $X$ as in Theorem~\ref{teo1} and assume that $\Sigma^c(X)=0$. Note that, by Robinson's version of Kupka-Smale theorem (\cite{R}) and well known results on linearizations, the hypotheses on the singularities of $X$ are generic. As $R$ admits a dominated splitting, a result of Vivier (Proposition 4.1 of~\cite{V}) implies that  $Sing(X)=\emptyset$, which is clearly an $C^1$-open condition. Now, a theorem of Zuppa (\cite{Z}) assures that $X$ can be $C^1$-approximated by a vector field $\tilde{X} \in \mathfrak{X}^{3}_{\mu}(M)$; moreover the hypotheses on $X$ guarantee that if $\tilde{X}$ is $C^1$-close to $X$ then $Sing(\tilde{X})=\emptyset$, $\tilde{X}$ is stably ergodic and $M$ is a partially hyperbolic set. If $\Sigma^c(\tilde{X})\not=0$ the proof ends. Otherwise we will construct a perturbation of $\tilde{X}$ to get a $C^2$-vector field $Y$ with $\Sigma^c(Y)\not=0$. This perturbation is made in two steps: first we obtain a kind of flowbox result in the setting of volume-preserving vector fields in order to trivialize coordinates (Lemma~\ref{cfbt}). Then, we make an explicit perturbation (Lemma~\ref{perturb}) and borrow Baraviera and Bonatti arguments (\cite{BB}) to show that this perturbation produces the desired vector field $Y$.

\begin{lemma}\label{cfbt}
Given a vector field $X\in\mathfrak{X}_{\mu}^{2}(M)$  and a
non-periodic point $p\in{M}$, there exists a conservative $C^{2}$
diffeomorphism $\Psi$ defined in a neighborhood of $p$ such that $T=\Psi_{*}X$, where $T= \frac{\partial}{\partial x_1}$.
\end{lemma}
\begin{proof}
Using Mosers's charts (see~\cite{Mo}, Lemma 2) we assume that $p=\vec{0}$ and that $X(p)=(1,0,...,0)$. For a small $r>0$ let $B_{r}(p)$ denote the $(n-1)$-dimensional ball centered at $p$ and of radius $r$ contained in $N_p$, and define the functions
$f,g:B_{r}(p)\rightarrow{\mathbb{R}}$ such that $f:=1$ and $g(x_2,...,x_n):=X_{1}(0,x_2,...,x_n)$, where $X_{1}$ is the first coordinate of the vector field $X$. Applying
Theorem 1 of~\cite{DM} we obtain a $C^{2}$ 
diffeomorphism
${\varphi}:B_{r}(p)\rightarrow{\varphi}(B_{r}(p))$ such that 
\begin{equation}\label{EDP}
g(\varphi(x_2,...,x_n))\text{det}D{\varphi}_{(x_2,...,x_n)}=\lambda
\end{equation}
for all $(x_2,...,x_n)\in{B_{r}(p)}$, where $\lambda=\int{g}/\int{f}$ and  $\varphi|_{\partial{B_{r}(p)}}=Id$.\\ Now we define the $C^{s}$ change of coordinates $\Phi:\mathbb{R}\times B_{r}(p)\rightarrow\mathbb{R}^{n}$ by $$\Phi((x_1,x_2,...,x_n))=X^{x_1\lambda^{-1}}(0,\varphi((x_2,...,x_n))).$$ 
By computing the jacobian of $\Phi$ at $(0,x_2,...,x_n)$ and using~(\ref{EDP}) we obtain:
\begin{equation}\label{origem}
\text{det}D{\Phi}_{(0,x_2,...,x_n)}=1, \forall(0,x_2,...,x_n)\in{\mathbb{R}\times\ B_{r}(p)}.
\end{equation}
Now we consider $(\overline{x}_{1},\overline{x}_2,...,\overline{x}_n)\in \mathbb{R}^+ \times B_{r}(p)$ and prove that $\text{det}D{\Phi}_{(\overline{x}_{1},\overline{x}_2,...,\overline{x}_n)}=1$.
It is clear that $\Phi(x_1,x_2,...,x_n)=X^{\overline{x}_{1}}[\Phi(x_1-\overline{x}_{1},x_2,...,x_n)]$ so $$D{\Phi}_{(x_1,x_2,...,x_n)}=DX^{\overline{x}_{1}}_{{\Phi}(x_1-\overline{x}_{1},x_2,...,x_n)}D{\Phi}_{(x_1-\overline{x}_{1},x_2,...,x_n)}.$$ 
Evaluated at $x_1=\overline{x}_{1}$ we obtain that 
$D{\Phi}_{(\overline{x}_{1},x_2,...,x_n)}=DX^{\overline{x}_{1}}_{{\Phi}(0,x_2,...,x_n)}D{\Phi}_{(0,x_2,...,x_n)}$. 
We use~(\ref{origem}) and the fact that the flow $X^{t}$ is volume-preserving to conclude that 
$\text{det}D{\Phi}_{(\overline{x}_{1},\overline{x}_2,...,\overline{x}_n)}=1$. Finally take $\Psi:=\Phi^{-1}$ and we obtain $T=\Psi_{*}X$.
\end{proof}

\medskip

Now we present the perturbation we mention before.

\begin{lemma}\label{perturb} Let $X\in\mathfrak{X}_{\mu}^{2}(M)$  and $\epsilon >0$. There exists $\xi >0$ with the following property: for any regular and non-periodic point $p$  and a two-dimensional vector space $V_p \subset N_p$, there is $r_0=r_0(p)>0$ such that for each $0<r\leq r_0$ exists $Y=Y_r\in\mathfrak{X}_{\mu}^{2}(M)$ which satisfies:
\begin{description}
\item (i) $Y$ is $\epsilon$-$C^1$-close to $X$;
\item (ii) $P_Y^1(p)_{|{W_p}}$ is the identity, where $W_p$ is the orthogonal complement of $V_p$ in $N_p$;
\item (iii) $P_Y^1(p)\cdot v=P_X^1(p) \circ R_{\xi} \cdot v$, $\forall v \in V_p$, where $R_{\xi}$ is the rotation of angle $\xi$ on $V_p$;
\item (iv) $X=Y$ outside the flowbox $X^{[0,1]}(B_{r}(p))$.
\end{description}
\end{lemma}
\begin{proof}
First, by using Lemma~\ref{cfbt} up to a change of coordinates $\Psi$, we can assume that $\Psi_ {*}X=T=\frac{\partial}{\partial x_1}$ in a neighborhood $\cal V$ of $p$. Actually, according to the proof of this Lemma, the neighborhood $\cal V$ can be chosen of the form $X^{[-1,1]}(B(p,r_0))$, for an appropriate small $r_0$ that depends on the point $p$. In these coordinates the space $N_p$ is the orthogonal space of  $\frac{\partial}{\partial x_1}$ and to simplify notation we still denote by $V_p$ and by $W_p$ the subspaces $D\Psi (p)(V_p)$ and $D\Psi (p)(V_p)$, respectively.

We fix two $C^{\infty}$ bump functions $\alpha,~\beta: \mathbb{R} \rightarrow [0,1]$ as follows. The function $\alpha$ is such that $\alpha(t)=0$ for $t\leq 0$, $\alpha(t)=1$ for $t\geq 1$, and $0\leq \alpha^{\prime}(t)\leq 2$, for all $t$; the function $\beta$  satisfies $\beta(t)=1$, for $|t|\leq \frac{1}{2}$, $\beta(t)=0$ for $t\geq 1$, and $|\beta^{\prime}(t)|\leq {4}$ for all $t$. For any $0<r<r_0 $ we define the map $\beta_r$ by $\beta_r(t)=\beta(\frac{1}{r}t)$.

For a fixed $r\in]0,r_0]$ we define $P: [-1,1]\times N_p \rightarrow N_p$ as follows. We fix an orthonormal basis of $V_p$, $\{u_2,u_3\}$, and, for $\theta \in [0, 2\pi]$, consider the rotation of angle $\theta$ whose matrix relative to this basis is   $$R_{\theta}=\begin{pmatrix} \cos(\theta)& -\sin(\theta)\\ \sin(\theta) & ~~\cos(\theta) \end{pmatrix}.$$ As $V_p \oplus W_p=N_p$, given $u\in N_p$ we write $u=u_V+u_W$, $u_V\in V_p$ and $u_W \in W_p$; for $t \in \mathbb{R}$ and $w \in N_p$ define 
$$P(t,w)(u)=R_{\alpha(t)\times \beta_r(\|w_V\|) \times \xi}(u_V)+u_W.$$
Let 
\begin{equation}\label{vari} 
\frac{\partial P(t,w)}{\partial t}=\left[\frac{\partial P(t,w)}{\partial t}\circ P(t,w)^{-1}\right]\cdot P(t,w)
\end{equation}
be the linear variational equation which will induce our perturbation. It is straightforward to see that
$$\frac{\partial P(t,w)}{\partial t}\circ P(t,w)^{-1}\cdot u=\alpha^{\prime}(t)\times \beta_r(\|u_V\|) \times \xi\begin{pmatrix} 0& -1\\ 1 & 0 \end{pmatrix}\cdot u_{V}.$$

We will define $Z$ on $\Psi({\cal V})$ such that 
\begin{equation}\label{Z}
DZ=\left[\frac{\partial P(t,w)}{\partial t}\circ P(t,w)^{-1}\right]. 
\end{equation}
For that we consider coordinates $(u_1,u_2,u_3,...,u_n)$ associated to an orthonormal basis of $\mathbb{R}^n$ such that $u_1$ is the  coordinate associated to the direction of the flow and $\{u_4,...,u_n\}$ is an orthonormal basis of $W_p$.
Now in order to get $Z$ satisfying ~(\ref{Z}) we take $$Z(u_1,u_2,u_3,...,u_n)=\alpha^{\prime}(u_1)\times \beta_r (\|(u_2,u_3) \|) \times \xi (0,-u_3,u_2,0,...,0).$$
Finally we define the vector field $Y$ by $Y=(\Psi^{-1})_*(T+Z)$ in $\cal V$, and by $Y=X$ in the complement of $\cal V$. We observe that $Y$ is a $C^2$ vector field (in fact $C^r$ if the initial $X$ is $C^r$), that condition \emph{(iv)} follows directly from the definition of $Y$  and note that a direct computation gives  $\text{div}(Y)=0$. 

To prove condition \emph{(i)} at each point $p$ we shrink, if necessary, the radius $r_0(p)<1$ in order to guarantee that $\frac{1}{2} \leq \frac{d_{\cal V}(X,W)}{d_{\Psi(\cal V)}(\Psi_*X, \Psi_*W)} \leq 2$, where $d_A(S,R)$ denotes the $C^1$ distance between the vector fields $S$ and $R$ over the set $A$. As $Y=(\Psi^{-1})_*(T+Z)$, from the expression of $Z$ it is not difficult to  prove that $d_{\Psi(\cal V)}(\Psi_*X, \Psi_*Y)$ is bounded by 
$ \text{max}\{ |\alpha^{\prime}| \beta_r  \xi,   |\alpha^{\prime \prime}| \beta_r  \xi,  
|\alpha^{\prime}| | \beta_r^{\prime}|  \xi r\}$. Since $| \beta_r^{\prime}| \leq 4/r$, $|\beta_r|\leq 1$ and $\text{max}|\beta_r^{\prime}|\geq 1$, from previous remark it follows that it is enough to choose 
$$ \xi \leq \frac{\epsilon}{2 \text{max}\{|\alpha^{\prime \prime}|,4|\alpha^{\prime}| | \beta_r^{\prime}|\}},$$ 
to assure that the $C^1$ distance between $X$ and $Y$ (which is equal to $d_{\cal V}(X,Y)$) is less than $\epsilon$.

To prove conditions \emph{(ii)} and \emph{(iii)} we observe that equation (\ref{vari}) and the definition of $Z$, see (\ref{Z}), imply that $P_{T+Z}^1(\Psi(p))=P(1,\vec{0})$; thus from the definition of $P(t,w)$ given above it follows that if $u \in V_p$ then $P_{T+Z}^1(\Psi (p))(u)=P(1,\vec{0})(u)=R_{\xi}(u)$, which implies \emph{(iii)} because $Y=(\Psi^{-1})_*(T+Z)$. Analogously we obtain \emph{(ii)} just by noting that $P_{T+Z}^1(\Psi (p))(u)=P(1,\vec{0})(u)=u$, if $u \in W_p$.
\end{proof}
\end{section}

\begin{section}{Proof of Theorem~\ref{teo1}}
\begin{subsection}{A key result and a simplified scenario}
Given a vector field $X$ satisfying the hypotheses of Theorem~\ref{teo1} using Lemma~\ref{perturb} we can exhibit a $C^2$-vector field, $Y$, arbitrarily $C^1$-close to $X$, without singularities, which is ergodic and partially hyperbolic. Moreover the two dimensional subspace $V_p$ (where the rotation acts) is chosen inside the central-unstable subbundle of $T_{p}M$, where $p$ is an Oseledets's point such that the orbit of $p$ is  non-periodic. Now our aim is to obtain the following result which ends the proof of the Theorem.

\begin{proposition}\label{main} Given $X$ and $p$ as above and an arbitrary small $r>0$ there exists $Y\in \mathfrak{X}_{\mu}^1(M)$ satisfying the conclusions of Lemma~\ref{perturb} and such that
$$\Sigma^c(Y) > \Sigma^c(X).$$
\end{proposition}

The proof of this result follows the strategy developed by Baraviera and Bonatti (\cite{BB}). Therefore here, for completeness and clearness, we do the proof in a very simplified scenario but give special attention to the construction of the linear differential system associated to the perturbation $Y$ of $X$. This system plays a crucial role in the estimates of the unstable Lyapunov exponents of $Y$. Concerning with the simplified setting we will restrict to the four dimensional case, we assume that the subbundles $N^u$, $N^c$ and $N^s$ are nontrivial (thus they are one-dimensional) and assume also that, in a neigbourhood of non-periodic point $p$, they are given by $\frac{\partial}{\partial x_2}$, $\frac{\partial}{\partial x_3}$ and $\frac{\partial}{\partial x_4}$, respectively.

\end{subsection}
\begin{subsection}{Construction of an useful linear differential system}
In the fixed setting the space $V_p$ is the space given by the directions $u_2$ (the unstable direction) and $u_3$ (the central direction), and $W_p$ is the space generated by $u_4$ (the stable direction).
We recall that $B_{r}(p)$ denotes the $(n-1)$-dimensional ball centered at $p$ and of radius $r$ contained in $N_p$.

For $q \in M$ let $u_2(q)$ denote a unit vector of $N^u_q$; write $$P_X^1(q)\cdot u_2(q)=\lambda(q)u_2(X^1(q)),$$ where $\lambda(q)=|\text{det}P_X^1(q)_{|N^u_q} |=J^u_X(q)$. The linear differential system is defined in the following way
$$\Phi(q,u_2(q))=(Y^1(q), \gamma(q) u_2(Y^1(q)),$$ where $\gamma(q)$ is equal to 
\begin{itemize}
\item[i)] $\lambda(q)$, if $q \notin X^{[-1,1]}(B_{r}(p))$,
\item[ii)] $\lambda (X^{-1}Y^1(q)) \cos(\beta_r(\|\overline{q}_V\|)\times \xi)$, if $q \in X^s(B_{r}(p))$, for some $s \in [0,1]$, where $\overline{q}=X^{-s}(q)$,
\item[iii)] $\lambda (q) \times A({q}^{*})$, if $q \in X^{s}(B_{r}(p))$, for some $s \in [-1,0[$, where ${q}^{*}=X^{-1-s}(q)$, where $A(\cdot)$ is defined bellow.
\end{itemize}
In order to define $A(q)$ let us fix $q \in X^{-1}(B_{r}(p))$ and define $\tau(q)$ as the least positive real number  such that there exists $\tilde{q} \in B_{r}(p)$ with $Y^{\tau (q)}(\tilde{q})=q$. Observe that $Y^{[\tau (q)]+1}(\tilde{q}) \in X^{[-1,0[}(B_{r}(p))$. If there is no recurrence, say if does not exist such a $\tau(q)$, then we will take $A(q)=1$. Take $\hat{q}=X^{-1}Y^1(\tilde{q})$; it follows that $q=X^{\tau (q)}(\hat{q})$. Moreover we also have $$P(1,\tilde{q})(u_{2}(\tilde{q}))=\cos(\beta_r(\|\tilde{q}_V\|)\times \xi) u_{2}(\hat{q})+v(\hat{q}),$$
where $v(\hat{q})=\sin(\beta_r(\|\tilde{q}_V\|)\times \xi) u_{3}(\hat{q})$
Therefore,
\begin{eqnarray*}
 P_{Y}^{\tau(q)}(\tilde{q})\cdot u_{2}(\tilde{q})&=&P_{Y}^{\tau(q)-1}(Y^{1}(\tilde{q}))\circ P_{Y}^{1}(\tilde{q})\cdot u_{2}(\tilde{q})=\\
&=& P_{Y}^{\tau(q)-1}(Y^{1}(\tilde{q}))\circ P_{X}^{1}(\tilde{q})\circ P(1,\tilde{q})(u_{2}(\tilde{q}))=\\
&=& P_{X}^{\tau(q)-1}(X^{1}(\tilde{q}))\circ P_{X}^{1}(\tilde{q})[\cos(\beta_r(\|\tilde{q}_V\|)\times \xi) u_{2}(\hat{q})+v(\hat{q})]=\\
&=& P_{X}^{\tau(q)}(\tilde{q})[\cos(\beta_r(\|\tilde{q}_V\|)\times \xi) u_{2}(\hat{q})+v(\hat{q})]=\\
&=& \cos(\beta_r(\|\tilde{q}_V\|) \xi) a(\hat{q})\left(\prod_{i=0}^{[\tau(q)]-1} \lambda(X^{i}(\hat{q}))\right) u_{2}(q)+P_{X}^{\tau(q)}(\tilde{q})\cdot v(\hat{q}),\\
\end{eqnarray*}
where $a(\hat{q})=|\text{det}P_X^{(\tau(q)-[\tau(q)])}{|N^u_{X^{[\tau]}(\hat{q})}} |$.
As $M$ is compact and the linear Poincar\'e flow is continuous it follows that there exists $C>0$ such that $C^{-1}<a(x)<C$, $\forall x \in M$. Moreover we also have that the linear differential system satisfies:
$$\Phi^{\tau(q)}(\tilde{q},u_{2}(\tilde{q}))=\left(Y^{\tau(q)}(\tilde{q}),\,\,\,\cos(\beta_r(\|\tilde{q}_V\|) \xi) a(\hat{q})\left[\prod_{i=0}^{[\tau(q)]-1} \lambda(X^{i}(\hat{q}))\right] u_{2}(q)\right).$$

We note that the abstract linear differential system differs from the linear Poincar\'e flow (associated to the perturbed flow $Y^t$) along the unstable direction by the projection in this direction of the vector $P_{X}^{\tau(q)}(\tilde{q})\cdot v(\hat{q})$. Actually we will take $A(q)$ as the ratio of the norms of these two vectors. So let us define
\begin{equation}\label{Aq}
A(q)=\frac{\| \Pi^u_{q} ( P_Y^{\tau(q)} (\tilde{q})\cdot u_2(\tilde{q}))\|}{\|\Phi^{\tau(q)}_2(\tilde{q}, u_2(\tilde{q}))\|}=1+\frac{\langle w_2(q),u_2(q) \rangle}{|\cos(\beta_r(\|\tilde{q}_V\|) \xi) a(\hat{q})\prod_{i=0}^{[\tau(q)]-1} \lambda(X^{i}(\hat{q}))|},
\end{equation}
where $\Pi^u_{q}$ is the projection from $N_ q$ on $N^u_q$ along the new central bundle, $\Phi_2$ denotes the second component of $\Phi$ and $w_2(q)=\Pi^u_{q}(P_{X}^{\tau(q)}(\tilde{q})\cdot v(\hat{q}))$. 

\bigskip

We observe that we constructed a measurable linear differential system over a Lebesgue invariant flow $Y^{t}$, hence by applying the Oseledets's Theorem (see for exam\-ple~\cite{JPS}) we conclude that for Lebesgue a.e. point the system has a Lyapunov exponent. From the Birkoff's Ergodic Theorem it follows that this Lyapunov exponent is equal to $\int_M \log(\gamma(q))d\mu(q)$.
\end{subsection}
\begin{subsection}{End of the proof}
Arguing exactly as in Lemma 1.4 of ~\cite{BB} we can prove that the Lyapunov exponent of the linear differential system  is equal to the unstable Lyapunov exponent of $Y^t$; hence, using   Birkoff's Ergodic Theorem again, we get 
\begin{equation}\label{le}
\int_M \log (J^u_Y(q)) d\mu(q)=\int_M \log(\gamma(q))d\mu(q).
\end{equation}

It is straightforward to see that
\begin{eqnarray*}
I(r)&=&\int_{B_{r}(p)} \log |\langle P(1,w)(\partial / \partial u_2), \partial / \partial u_2 \rangle|d\hat{\mu}(w) =\\
&=&\int_{B_{r}(p)} \log|\cos(\beta_r(\|w_V\|)\times \xi)|d\hat{\mu}(w)=\\
&=&\text{vol}(B_{r}(p))\int_{B_{1}(p)} \log|\cos(\beta(\|w_V\|)\times \xi)|d\hat{\mu}(w)=\text{vol}(B_{r}(p)) I(1)<0.
\end{eqnarray*}

Define $\Sigma_X^u=\int_M \log (J^u_X(q)) d\mu(q)$ and $\Sigma_\Phi^u=\int_M \log (\gamma(q)) d\mu(q)$. Let us estimate $\Sigma_X^u-\Sigma_\Phi^u$. First note that these two numbers coincide outside the flowbox  $X^{[-1,1]}(B_r(p))$. So let us begin by computing its difference on $X^{[0,1]}(B_r(p))$. Let $\hat{\mu}$ be the three dimensional  ``Lebesgue measure'' associated to $\mu$ and defined on normal sections. Using that $X^{-1}Y^1$ is a volume preserving diffeomorphism which leaves  $X^{[0,1]}(B_r(p))$ invariant it follows that
\begin{eqnarray*}
&&~~~\int_{X^{[0,1]}(B_r(p))}\log (J^u_X(q))- \log (\gamma(q)) d\mu(q)=\\
&=&-\int_{B_r(p)}\int_0^1 \log  (\cos(\beta_r(\|X^{-s}(q_V)\|) \xi))dsd\hat{\mu}(\overline{q})=\\
&=&-\int_{B_r(p)} \log  (\cos(\beta_r(\|\overline{q}_V\|) \xi))d\hat{\mu}(\overline{q})=\\
&=&-\text{vol}(B_{r}(p)) I(1).
\end{eqnarray*}
On the other hand we have
\begin{eqnarray*}
\int_{X^{[-1,0[}(B_r(p))}\log (J^u_X(q))- \log (\gamma(q)) d\mu(q)&=&-\int_{X^{[-1,0[}(B_r(p))}\log A(q^{*})d\mu(q)=\\
&=&-\int_{X^{-1}(B_r(p))}\log A(q^{*})d\hat{\mu}(q^{*}).
\end{eqnarray*}
Therefore we obtain
\begin{eqnarray*}
\Sigma_X^u-\Sigma_\Phi^u &\geq &-\text{vol}(B_{r}(p)) I(1) - \text{vol}({X^{-1}(B_r(p))}) \times \underset{q^{*} \in X^{-1}(B_r(p))}{\text{max}} \log(A(q^{*}))=\\ 
&=&-\text{vol}(B_{r}(p))\left(I(1) + \frac{\text{vol}({X^{-1}(B_r(p))})}{\text{vol}(B_{r}(p))} \underset{q^{*} \in X^{-1}(B_r(p))}{\text{max}} \log(A(q^{*})) \right).
\end{eqnarray*}
As $\frac{\text{vol}({X^{-1}(B_r(p))})}{\text{vol}(B_{r}(p))}$ is uniformly bounded, for $p \in M$ and $r \leq r_0$, and $I(1)<0$, to prove that $\Sigma_X^u-\Sigma_\Phi^u >0$  it is enough to show that $\underset{q^{*} \in X^{-1}(B_r(p))}{\text{max}} \log(A(q^{*}))$ converges to zero as $r$ tends to zero. For that it  suffices to prove  that there exist $\sigma \in ]0,1[$ and $C>0$ such that for any small $r$ and $q^{*} \in X^{-1}(B_r(p))$ one has $|A(q^{*})-1| \leq C\sigma^{\tau_r}$, where $\tau_r$  is the smallest return time from $B_r(p)$ to $B_r(p)$. Now, as $\langle w_2(q),u_2(q) \rangle \leq \|w_2(q)\|$,  $|\cos(\beta_r(\|\tilde{q}_V\|) \xi)| \approx 1$ and $a(\hat{q})$ is bounded away from zero, from ~(\ref{Aq}), we get that there exists $C_1>0$ such that 
$$|A(q^{*})-1|\leq \frac{C_1  \|w_2(q)\|}{\prod_{i=0}^{[\tau(q)]-1} \lambda(X^{i}(\hat{q}))}.$$
In the next step we use the domination structure:
$$\| P_X^t(\hat{q}) \cdot (w_c) \| \leq  C_2\sigma^t \| P_X^t(\hat{q}) \cdot (w_u)\| \leq C_3 \sigma^t \prod_{i=0}^{[\tau(q)]-1} \lambda(X^{i}(\hat{q})),$$
for some constant $C_3>0$, where the unit vector $w_c$, respectively $w_u$, belongs to $N^c(\hat{q})$, respectively  $N^u(\hat{q})$.
Since $w_2(q)=\Pi^u_{q}(P_{X}^{\tau(q)}(\tilde{q})\cdot v(\hat{q}))$ and $\tau(q) \geq \tau_r$ we obtain
$$|A(q^{*})-1|\leq \frac{C_1 C_3 \sigma^{\tau(q)}\| \Pi^u_{q}(P_{X}^{\tau(q)}(\tilde{q})\cdot v(\hat{q}))\|}{\|P_{X}^{\tau(q)}(\tilde{q})\cdot v(\hat{q})\|} \leq C_1 C_3 \sigma^{\tau_r}. $$
Since $\tau_r$ tends to infinity as  $r$ goes to zero it follows that if $r$ is small then $\Sigma_X^u-\Sigma_\Phi^u >0$, where, we recall, $\Phi$ is the linear differential system associated to $Y$ and this vector field is the perturbed one associated to $r$.
To end the proof of the proposition we observe that $\Sigma_Y^u=\Sigma_\Phi^u$; moreover, by the volume preserving assumption, 
$$\Sigma^u_X+\Sigma^c_X+\Sigma^s_X=0=\Sigma^u_{Y}+\Sigma^c_{Y}+\Sigma^s_{Y}.$$
As in this simplified scenario the negative Lyapunov exponents are unchanged, that is $\Sigma_Y^s=\Sigma_X^s$, we conclude that $\Sigma_Y^c>\Sigma_X^c$, where $\Sigma_X^c$ is the sum with multiplicities of the central Lyapunov exponents.
\end{subsection}
\end{section}
\bigskip

\textbf{Acknowledgements:} We would like to thank Alexandre Baraviera for suggesting us the central question we treated in this work.

\flushleft
{\bf M\'ario Bessa} \ \  (bessa@impa.br)\\
CMUP, Rua do Campo Alegre, 687 \\ 4169-007 Porto \\ Portugal\\

\medskip

\flushleft
{\bf Jorge Rocha}  \ \  (jrocha@fc.up.pt)\\
DMP-FCUP, Rua do Campo Alegre, 687 \\ 4169-007 Porto \\ Portugal\\

\end{document}